\documentclass[10pt]{article}
\usepackage{amssymb}
\usepackage{latexsym}
\usepackage{euscript}

\newcommand{\pp}{\mathbb{P}}

\newcommand{\C}{\mathbb{C}}
\newcommand{\R}{\mathbb{R}}

\newcommand{\Z}{\mathbb{Z}}
\newcommand{\sss}{\mathbb{S}}
\newcommand{\T}{\mathbb{T}}

\DeclareFontFamily{OT1}{rsfs}{}
\DeclareFontShape{OT1}{rsfs}{n}{it}{<-> rsfs10}{}
\DeclareMathAlphabet{\curly}{OT1}{rsfs}{n}{it}

\newcommand{\Sym}{\mathrm{Sym}}

\newcommand{\M}{\mathcal{M}}

\newcounter{Universal}[section]
\renewcommand{\theUniversal}{\thesection.\arabic{Universal}}

\newenvironment{Plain}{\refstepcounter{Universal} \par \vspace{0.5cm}
\noindent {\bf (\theUniversal)}\ }{\par \vspace{0.5cm}}
\newenvironment{Italic}{\refstepcounter{Universal} \par \vspace{0.5cm}
\noindent {\bf (\theUniversal)}\ \it}{\par \vspace{0.5cm}}  
  
\newenvironment{Thm}{\begin{Italic}{\sc Theorem: }}{\end{Italic}}
\newenvironment{Prop}{\begin{Italic}{\sc Proposition:} }{\end{Italic}}
\newenvironment{Defn}{\begin{Italic}{\sc Definition: } }{\end{Italic}}  
\newenvironment{Cor}{\begin{Italic}{\sc Corollary: }}{\end{Italic}}
\newenvironment{Lem}{\begin{Italic}{\sc Lemma: }}{\end{Italic}}
\newenvironment{Conj}{\begin{Italic}{\sc Conjecture: }}{\end{Italic}}
\newenvironment{Pf}{\par \noindent{\sc Proof:} }{\quad $\blacksquare$ \par
\vspace{0.5cm}}

\newcounter{enum}

\newenvironment{Eqn}{\refstepcounter{Universal} $$} {\eqno \mathrm{
(\theUniversal)} $$} 

\title{On moduli spaces of symplectic forms}
\author{Ivan Smith}
\date{New College, Oxford}

\begin{document}
\maketitle
\thispagestyle{empty}

\begin{abstract}
\noindent  We prove that there are simply-connected 
four-manifolds which admit $n$-tuples of symplectic forms whose
first Chern classes have pairwise different divisibilities in integral
cohomology.  It follows that the moduli spaces of symplectic
forms modulo diffeomorphisms on the manifolds are disconnected.
\end{abstract}

%%%%%%%%%%%%%%%%%%%%%%%%%%%%%%%%%%%%%%%%%%%%%%%%%%%%%%%%%%%%%%%%%%%%

\section{Introduction}

This paper provides examples of symplectic manifolds with disconnected
moduli spaces of symplectic forms.  Recall that
two symplectic forms $\omega$ and $\omega'$ on a symplectic manifold
  $X$ are \emph{deformation 
  equivalent} if there is some path of symplectic forms interpolating
between them.  We allow the cohomology class to vary both along the
path and at the endpoints.  The forms $\omega$ and $\omega'$ are
\emph{diffeomorphism related} if there is some diffeomorphism $\phi: X
\rightarrow X$ for which $\phi^* \omega = \omega'$; we do not insist
that $\phi$ is isotopic to the identity.  We write $\M(X)$ to
denote the moduli space of symplectic forms on $X$ up to
diffeomorphism.

\begin{Defn} \label{equivalent}
  The forms $\omega$ and $\omega'$ are \emph{equivalent} symplectic
  forms on $X$ if they lie in the same equivalence class under the
  equivalence relation generated by deformations and diffeomorphisms.
\end{Defn}

\noindent  Thus equivalence classes of symplectic form are
indexed by $\pi_0 (\M (X))$.
This is a suitable notion of equivalence for forms when we are
interested in essentially topological properties of the underlying
symplectic manifolds.  Note that even coarse symplectic invariants of
$X$, such as total volume, are not preserved by an equivalence of
forms.  

\vspace{0.2cm}

\noindent In dimensions $4k$ with $k>1$ there are many examples of
manifolds with inequivalent symplectic forms, which can be
distinguished by the divisibility of the first Chern class in integral
cohomology.  (In \cite{Ruan}, using Gromov-Witten invariants, Yong-Bin
Ruan gave sophisticated examples
of inequivalent forms -  with the same Chern classes -  on
algebraic $3$-folds; notice that he uses the term
``deformation equivalent'' for our ``equivalent''.  Ruan also 
gave examples \cite{Ruan2} of diffeomorphisms of $3$-folds which do
intertwine Chern 
classes but which do not induce an equivalence of forms.)
To see this, note that given homeomorphic four-manifolds $X$ and $Y$,
the products $X \times X \ = \ Y \times Y \ = \ Z$ are
diffeomorphic.  Suppose 
$X$ is a minimal complex surface of general
type with first Chern class divisible by any odd $n>1$, and $Y$ is a
non-minimal symplectic
four-manifold homeomorphic to $X$.  (It is straightforward to construct
such pairs using Freedman's theorems.)  The product symplectic form
$\omega_X \oplus \omega_X$ on $Z$ has odd $c_1$ whereas the form
$\omega_Y \oplus \omega_Y$ has $c_1$ prime.  Taking additional
products with further copies of the manifolds extends this to all higher
dimensions $4k$.  However, in dimension four until recently no manifolds
with inequivalent symplectic forms were known.  This question was
settled by McMullen and Taubes in \cite{McM-Taubes}, but the forms
were distinguished by gauge theory techniques specific to four
dimensions.  In this 
paper we provide four-dimensional examples on a similar footing to the
higher dimensional examples sketched above.  Our main
result is the following

\begin{Thm} \label{alldimensions}
For each $n \geq 2$ there is a simply-connected manifold $X^4$ which
admits $n$ distinct
symplectic forms whose first Chern classes are of pairwise different
divisibilities in $H^2 (X; \Z)$.  In particular for such an $X$ we have
$|\pi_0 (\M (X))| \geq n$.  
\end{Thm}

\noindent We note that a recent result of Li and Liu \cite{Li-Liu}
shows that on any symplectic four-manifold, there are at most finitely many
orbits of the diffeomorphism group on the set of first Chern classes
of symplectic structures.  In some sense the above result is therefore
optimal\footnote{After this paper was accepted, the author learned
  that independently Stefano Vidussi (\emph{Homotopy $K3$'s with several
    symplectic structures}) has used gauge theory methods to show
  the existence, for any $n$, of homotopy $K3$ surfaces with $n$ inequivalent
  symplectic structures.}. (In contrast to the work of Ruan mentioned
above, the 
existence of inequivalent symplectic forms not
distinguished by the first Chern class remains open in four
dimensions.)  Theorem (\ref{alldimensions}) perhaps casts a new light, via
taking products, on a
conjecture of Donaldson (\ref{conjecture} and \ref{conjecturefails}).
The proof shall require only elementary arguments.  Before
giving more details, we set the result in context. 

\vspace{0.2cm}

\noindent $\bullet$ \textbf{Moduli spaces of symplectic forms:}
Following work of Taubes relating Seiberg-Witten invariants to
pseudoholomorphic curves, there are a few symplectic
four-manifolds for which the moduli space of symplectic forms is
entirely understood.  In all such cases, the understanding comes from
a uniqueness theorem:  for instance, every symplectic form on $\C
\pp^2$ is standard up to scale and diffeomorphism \cite{Taubes}.
Lalonde and McDuff have proved similar 
results for other rational and ruled manifolds \cite{Lal-McD}.  Moreover in
four dimensions many symplectic invariants are largely determined by
the smooth structure of the underlying manifold, most famously for certain
Gromov invariants counting pseudoholomorphic curves \cite{Taubes}.  Recently
McMullen and Taubes \cite{McM-Taubes} gave the first examples of symplectic
four-manifolds for which there are symplectic forms which are not
equivalent to each other under any composition of deformations and
pullbacks by diffeomorphisms.  Their examples were distinguished by a
classical invariant of a symplectic structure, the first Chern class,
but via non-classical methods.  First Chern classes of symplectic
structures on manifolds with $b_+ > 1$ are always Seiberg-Witten basic
classes.  Diffeomorphisms of a manifold preserve the set of basic
classes inside $\R^{b_2} = H^2 (X; \R)$;  since diffeomorphisms act on
this space linearly, they preserve the polytope which is the convex
hull of the basic classes.  It follows that they must preserve the
sets of vertices of this polytope of any fixed valence;  McMullen and
Taubes find two symplectic structures on a fixed manifold for which
the first Chern classes represent vertices of distinct valence.  

Detecting the non-transitivity of diffeomorphisms on basic classes is
of great interest in its own right, but it is natural to hope that one
might distinguish symplectic structures on simply connected manifolds
using purely classical methods.  When a manifold contains a Lagrangian
torus of square zero which is homologically essential, there are
perturbations of the symplectic form making the torus a symplectic
submanifold with \emph{either choice of orientation}.  This
flexibility, which underlies the McMullen-Taubes construction in
different language, will enable us to give new and easy examples of
simply connected four-manifolds for which the moduli space of symplectic
forms is disconnected.  (Embedded surfaces of square zero were also
used in \cite{ivansympsub} to give a topological construction of
exotic symplectic submanifolds on complex surfaces.)

\vspace{0.2cm}

\noindent $\bullet$ \textbf{Donaldson's Conjecture:}
According to (\cite{McD-S}; p.437) Simon Donaldson
has formulated the following remarkable conjecture relating
diffeomorphisms of symplectic four-manifolds to equivalences of
symplectic six-manifolds.  Write $\omega_{st}$ for the standard
K\"ahler form on projective space.

\begin{Conj} \label{conjecture}
Let $(X, \omega_X)$ and $(Y, \omega_Y)$ be symplectic four-manifolds
which are homeomorphic.  Then they are diffeomorphic if and only if
the product symplectic forms $\omega_X \oplus \omega_{st}$ and
$\omega_Y \oplus \omega_{st}$ are equivalent on $X \times \sss^2
\cong_{\mathrm{diff}} Y \times \sss^2$.
\end{Conj}

\noindent We stress that we have nothing to say about this directly;
however,  it 
follows from our examples that the sphere $\sss^2$ is playing a special
role here which is not immediately apparent from the formulation of
(\ref{conjecture}).  From (\ref{alldimensions}) we see

\begin{Cor} \label{conjecturefails}
If in the statement of Donaldson's conjecture one replaces $\sss^2$ by 
a symplectic manifold with $c_1$ divisible by $n > 2$, then the
analogous conjecture fails.
\end{Cor}

\noindent Hence we cannot replace $\sss^2$ by the torus $\T^2$ or the
$K3$ surface.  A more vivid statement is provided by the following
contrast.  Donaldson's conjecture would imply that
the obvious symplectic structures on $X \times (\pp^1 \times \pp^1) \
= \ Y \times (\pp^1 \times \pp^1)$ are equivalent, whereas
(\ref{conjecturefails}) shows that if we 
replace $\pp^1 \times \pp^1$ with $\Sym^2 (\pp^1) \cong \C \pp^2$ then
this need not be true.

\vspace{0.2cm}

\noindent The heart of (\ref{alldimensions}) is the following
sufficient condition for a four-manifold to admit inequivalent
symplectic forms.  
Write $E(1)$ for the rational elliptic surface \cite{GompfS}.

\begin{Thm} \label{SLtorimflds}
Let $X$ be a symplectic four-manifold such that

\begin{enumerate}
\item there are square-zero tori $T_1, \ldots, T_r \subset X$ which span an
  $r$-dimensional subspace in $H_2 (X; \R)$, for some $r>1$;
\item $T_1$ is represented by a Lagrangian submanifold,
  $T_{i>1}$ and $\sum_i T_i$ are represented by symplectic
  submanifolds (tori), all disjointly embedded;
\item the first Chern class $c_1 (X)$ is divisible by $n$ for some
  $n>2$, whilst\footnote{We use the symbol $a \not\mid b$ to show that
    $a$ does not divide $b$.} $n \not\mid 2d$ where $d$ is the
  divisibility of $[T_1] \in H_2 (X; \Z)$.
\end{enumerate}

\noindent Then the manifold $Z$ given by fibre summing $X$ with a copy of
$E(1)$ along each of $T_1, \ldots, T_r$ and along $(n-1)$ parallel
copies of $T_1 + \cdots +T_r$ admits inequivalent symplectic
structures.  
\end{Thm}

\noindent The smooth structure on $Z$ depends on a choice, once and
for all, of trivialisations of the normal bundles of the tori.  If
$\pi_1 (X)$ is normally generated by the images of the $\pi_1 (T_i)$
then $Z$ is simply connected.  Here is an example.  Let $\T^4$ describe the
usual Euclidean four-torus $\R^4 / \Z^4$.  With co-ordinates $x,y,z,t$
on $\R^4$ there are two-tori in $\T^4$ which come from projecting the
planes described by the pairs of lines

$$\langle x,t \rangle; \ \langle y,t \rangle; \ \langle z,t \rangle; \ 
\langle x=y=z, t \rangle.$$

\noindent Call these respectively $T_x, T_y, T_z, T_w$.  All the tori are
oriented by taking the usual orientations on the defining lines in
$\R^4$ (pointing out towards $+ \infty$).  We can perturb the tori to
obtain five disjoint oriented tori inside $\T^4$; two parallel copies
of $T_w$ and one of each of $T_x, T_y, T_z$.  Let $Z$ denote the
simply connected smooth four-manifold obtained by forming the normal
fibre sum of $\T^4$ with five copies of $E(1)$, gluing a standard
complex fibre $F \subset E(1)$ to each of the five tori above.

\begin{Cor} \label{twoformsdiffer}
  The simply-connected manifold $Z$ admits two symplectic structures
  which are not equivalent under any sequence of deformations or
  diffeomorphisms.  Indeed there are symplectic structures $\omega_1$
  and $\omega_2$ on $Z$ for which $c_1 (TZ; \omega_1)$ is divisible by
  $3$ in integral cohomology but $c_1 (TZ; \omega_2)$ is prime.
\end{Cor}

\noindent Since the construction allows for some variation, infinitely many
homeomorphism types of simply connected four-manifold displaying a similar
phenomenon can be obtained.  To obtain a manifold with $N$
inequivalent symplectic forms, take $n=p_1 \cdots p_N$ to be a product of $N$
distinct odd primes in (\ref{SLtorimflds}).  Then by varying the number of
fibre sums performed with the ``wrong'' orientation along the
Lagrangian torus, one can obtain symplectic forms on the final
manifold $Z$ for which the first Chern class is divisible by precisely
one of the primes $p_i$.  Thus (\ref{alldimensions}) will be a
direct consequence of the above results.  We will address the proofs of
(\ref{SLtorimflds} and \ref{twoformsdiffer}) in the subsequent
sections.  Using these, and the remarks at the end of the paper, it is
straightforward to deduce
(\ref{conjecturefails}); we leave the details to the reader.

\vspace{0.3cm}

\noindent \textbf{Acknowledgements:}
The McMullen-Taubes construction was cast
in the light presented here during a lecture by Ron Fintushel, and I
am grateful to him for helpful correspondence.  I am also indebted to Bob
Gompf, who offered valuable comments on an earlier draft of the
paper. 

%%%%%%%%%%%%%%%%%%%%%%%%%%%%%%%%%%%%%%%%%%%%%%%%%%%%%%%%%%%%%%%%%%%%%%

\section{Lagrangian fibre sums}

The fibre-sum construction for symplectic manifolds was pioneered by
Gompf in \cite{Gompf} who used it to construct symplectic manifolds
displaying many new phenomena.  Gompf also noticed that given a
homologically essential Lagrangian submanifold $L$ of a symplectic
four-manifold $X$, one could perturb the ambient symplectic form
$\omega$ by an arbitrarily small amount so as to make the form
non-vanishing on $L$.  In the presence of a fixed orientation on $L$,
there is some choice here; one can perturb $\omega$ to be either
positive or negative on $L$.  Recently this construction, viewed from
a somewhat different perspective, has been used by McMullen and Taubes
to produce symplectic manifolds with inequivalent symplectic forms.
Their examples were distinguished using Seiberg-Witten theory, and the
action of a diffeomorphism of a four-manifold on the convex hull of
the Seiberg-Witten basic classes.  Via the description in terms of
fibre sums along Lagrangians, one can obtain simpler examples of the
same phenomena; we present one class here.  The author is very
grateful to Ron Fintushel for a lucid explanation of the original work
of McMullen and Taubes.

\vspace{0.2cm}

\noindent We begin with some more general remarks on the fibre sum operation in
our context.  Let $(X,T)$ and $(Y,T')$ be symplectic pairs, where
$X,Y$ are symplectic four-manifolds and $T,T'$ are embedded tori of
square zero on which the symplectic forms are non-degenerate.  The
symplectic structures define orientations on all four manifolds, and
hence orientations on the normal bundles of $T$ in $X$ and $T'$ in
$Y$.  Then given any orientation-preserving diffeomorphism $f$ of $T$
and $T'$, we may lift to an orientation-reversing diffeomorphism of
the normal bundles and use this to glue the manifolds $X$ and $Y$
along a neighbourhood of the embedded surfaces.  Explicitly, choose a
marking of $T$ and $T'$ with a standard torus so the diffeomorphism
$f$ becomes the identity.  Then we form

$$X \sharp_{(T:T')} Y \ = \ X \backslash \nu_X T \ \cup_{\tau} \ Y
\backslash \nu_Y T'$$

\noindent where $\tau: T \times D^2 \rightarrow T' \times D^2$ is the
map taking $(t,z) \mapsto (t, \overline{z})$.  Since complex
conjugation reverses orientation on the disc and on its boundary
circle, there is a natural orientation induced on the final manifold.
The diffeomorphism type of the final manifold depends \emph{a priori}
on the diffeomorphism $f$ and the choice of lift to an identification
of the normal bundles.  In our applications, the normal bundles of the
tori will be canonically trivial.  If in addition we suppose that $(Y,
T')$ is an elliptic surface $(E(1), F)$ with an embedded complex
fibre, then in fact all choices of $f$ give orientation-preserving
diffeomorphic manifolds.  This is because of the standard

\begin{Lem} \label{extend}
  Let $E(1)$ denote the complex elliptic surface given by blowing up
  the nine base-points of a generic pencil of cubics in $\C \pp^2$.
  Let $F$ denote a generic smooth complex fibre of the resulting
  elliptic fibration.  Then every orientation preserving
  diffeomorphism of $\nu_{E(1)}F$ may be extended to an orientation
  preserving diffeomorphism of $E(1)$.
\end{Lem}

\noindent Note importantly that in this case, \emph{the smooth
  structure on the fibre sum is determined completely by the
  orientation on $T$ and by the trivialisations of the normal
  bundles}.  The orientation of the torus in $E(1)$ is not important
since there are diffeomorphisms of the boundary of a neighbourhood of
the fibre which are orientation preserving but reverse the orientation
on the fibre.  Once all the choices for the smooth structure on the
sum are fixed, we can consider how to patch forms.  That is, the
symplectic structures themselves enter only in putting a symplectic
form on the sum.

We will make regular use of the following:

\begin{Prop} \label{firstchern}
  Let $W = X \sharp_{(T:T')} Y$ denote the symplectic sum of manifolds
  $X,Y$ along embedded symplectic tori of square zero.   Then
  
  $$c_1 (W) \ = \ c_1 (X) + c_1 (Y) - 2 \mathcal{P}([T=T']).$$

\noindent Here $\mathcal{P}(\cdot)$ denotes the Poincar\'e dual of a homology
class (we will suppress this from the notation henceforth).
\end{Prop}

\begin{Pf}
  The notation needs explanation, since $c_1 (X)$ and $c_1 (Y)$ are
  not naturally elements of $H^2 (W)$.  One can make sense of them in
  various ways.  Suppose inside 
  $X$ and $Y$ we have embedded  surfaces  which are 
  representatives for $c_1$ disjoint from the tori $T$ and $T'$.  In
  this case these surfaces define homology classes
  in the manifold $W$ and we take the Poincar\'e duals of these.
  The tangent bundle to an oriented torus is (canonically)
  trivial, so the
  canonical bundles of $X$ and $Y$ are trivial over neighbourhoods
  of $T$ and $T'$.  Hence we can obtain the surfaces we require by
  choosing smooth sections of $K_X^{-1}$ and $K_Y^{-1}$ which are constant and
  non-vanishing near $T$ and $T'$.  

\vspace{0.2cm}

\noindent There is apparently some choice. The trivialisations of the canonical
  bundles near the tori depend on a choice of trivialisation of the
  normal bundles of the tori.  For our later examples there are
  natural choices: the fibre in an elliptic fibration has a canonically
  trivial normal bundle, whilst any Lagrangian torus has a canonically
  trivial normal bundle since a (connected) choice of almost complex
  structure defines an isomorphism from the normal bundle to the
  tangent bundle.  Being more careful, recall that

$$H_2 (X \backslash \nu_X (T)) \ = \ H^2 (X \backslash \nu_X (T); \
\partial);$$ 

\noindent one can see that a choice of trivialisation of the normal
bundle to $T$ defines a relative Chern class in the group on the
right.  Moreover, we can use the given trivialisation of $\nu_X (T)$
and diffeomorphism $T \rightarrow T'$ 
to define one for $\nu_Y (T')$, hence a second relative Chern
class in $H^2 (Y \backslash \nu_Y (T'); \ \partial)$. Perform the
gluing relative to
these normal trivialisations; changing the
trivialisation affects the two relative Chern classes in cancelling
fashion.  Thus the final result is independent of choices.

\vspace{0.2cm}
 
\noindent  Interpreting the statement of the Proposition
appropriately, a proof runs 
  as follows.  We can choose smooth sections of the anticanonical bundles
  of each of $X, Y$ which are non-vanishing over a neighbourhood of
  the gluing tori.  Then the final Chern class differs from the sum
  by a contribution from the normal directions to the tori.  The
  tangent bundles to $X$ and $Y$ split as a product here, and the
  result follows from the usual formula for connect summing
  two-dimensional surfaces: 

$$e(\Sigma_1 \sharp \Sigma_2) \ = \ e(\Sigma_1) + e(\Sigma_2) - 2.$$

\noindent Equivalently, the result follows since complex
  conjugation on the 
  annulus has algebraically two negatively oriented fixed points;
  compare to a favourite K\"ahler example.
\end{Pf}

\noindent Now suppose $(X,T)$ is a pair comprising a symplectic four-manifold
and an \emph{oriented} Lagrangian torus $T$, which we suppose is
non-trivial in homology.  Gompf \cite{Gompf} notes that we can perturb the
symplectic form $\omega_X$ (by an arbitrarily small amount, and in
particular preserving the first Chern class of the symplectic
structure) to be non-degenerate on $T$.  Fixing a pair of oriented
tangent vectors $v_1, v_2$ to $T$ at a point $p \in T$, we then have
two possibilities: $\omega_X (v_1, v_2)_p >0$ or $\omega_X (v_1,
v_2)_p < 0$.  That is, the symplectic orientation induced on $T$ can
agree or disagree with the given fixed orientation.  In the first
case, we can form the symplectic sum of $(X,T)$ with $(E(1), F)$ in
the usual way, and we will find that (writing $T_+$ to note the
positivity)

\begin{Eqn} \label{choiceone}
c_1 (X \sharp_{(T_+:F)} E(1)) \ = \ c_1 (X) + c_1 (E(1); \omega_0) -
2[T_+=F]
\end{Eqn}

\noindent comparing with (\ref{firstchern}).  Here $\omega_0$ denotes the usual
K\"ahler form on the rational elliptic surface. However, for a
perturbation of $\omega_X$ for which $T$ is negatively oriented, if we
form the fibre sum of $X$ along $T$ with the \emph{given} and not
symplectic orientation for the normal bundle $\nu_{T/X}$, then we can
only sum with a pair $(Y,T')$ for which
$$\int_{T'} \omega_Y \ = \ \int_T \tilde{\omega}_X \ < 0$$
where the
$\tilde{\cdot}$ denotes the perturbed form.  Taking $Y = E(1)$ and $T'
= F$ we find that we can now form the symplectic sum if we take
$-\omega_0$, the negative of the usual K\"ahler form, as the
symplectic form on $E(1)$. Then looking at Chern classes we find

\begin{Eqn} \label{choicetwo}
c_1 (X \sharp_{(T_-:F)} E(1)) \ = \ c_1 (X) + c_1 (E(1); -\omega_0) -
2 [-T_-=-F].
\end{Eqn}

\noindent Since we are fibre summing along surfaces which have the
reversed orientation to the usual symplectic setting
(\ref{firstchern}), in the last term we reverse sign.  The following
is standard:

\begin{Lem} \label{changesign}
  For the rational elliptic surface we have
  
  $$c_1 (E(1); \omega_0) \ = \ [F], \qquad c_1 (E(1); -\omega_0) \ = \ 
  -[F]$$

\noindent identifying homology and cohomology via Poincar\'e duality as usual.
\end{Lem}

\noindent With these various comments in place, the theorem
(\ref{SLtorimflds}) is straightforward.

\begin{Thm}
  Suppose $X$ is a symplectic manifold as in the statement of
  (\ref{SLtorimflds}).  Pick an orientation on $T_1$ and form the
  oriented fibre sum of $X$ with $r+n-1$ copies of $E(1)$ as
  described.  Then the resulting manifold $Z$ admits inequivalent
  symplectic structures.  Moreover it is simply connected if $\pi_1
  (X)$ is normally generated by the images of the $\pi_1 (T_i)$.
\end{Thm}

\begin{Pf}
  We put two symplectic structures on $Z$, the first by perturbing
  $T_1$ to have the positive symplectic orientation and summing with
  $(E(1), \omega_0)$ and the second by perturbing in the other
  direction and summing with $(E(1), -\omega_0)$.  From the above
  lemmata we find that the two symplectic structures on $Z$ have Chern
  classes
  
  $$c_1 (Z)_+ \ = \ c_1 (X)-[T_1] -[T_2] - \cdots -[T_r] - (n-1) [T_1
  + \cdots + T_r]$$

\noindent in the positive case, and

$$c_1 (Z)_- \ = \ c_1 (X) + [T_1] - [T_2] - \cdots -[T_r] - (n-1) [T_1
+ \cdots + T_r]$$

\noindent in the negative case.  Since $n | c_1 (X)$ clearly $n | c_1
(Z; \omega_+)$; however, $c_1 (Z; \omega_-) = nA + 2[T_1]$ for $A =
(c_1 (X) / n) - \sum[T_i]$.  Recalling that $n \not\mid 2d$ for $d$
the divisibility of $[T_1]$, the Chern class $c_1 (Z)_-$ cannot be
divisible by $n$, and the result follows.  The statement on
fundamental groups is straightforward since $\pi_1 (E(1) \backslash F)
= 0$.
\end{Pf}

\noindent In the last section we provide some explicit examples.  For
simplicity, note that the divisibility condition $n \not\mid 2d$ will
always hold if $[T_1]$ is prime in homology, whilst $n | c_1 (X)$ is
trivial for any $n$ if $c_1 (X) = 0$.

%%%%%%%%%%%%%%%%%%%%%%%%%%%%%%%%%%%%%%%%%%%%%%%%%%%%%%%%%%%%%%%%%%%%%%%%%%%

\section{Inequivalent symplectic forms}

We develop the examples described in (\ref{twoformsdiffer}).  The
reader may wish to compare to the original construction of McMullen
and Taubes from \cite{McM-Taubes}.  Take co-ordinates $x,y,z,t$ on
$\R^4$ and write $\T^4 = \R^4 / \Z^4$.  Define the oriented tori $T_x,
T_y, T_z$ and $T_w$ as follows:

\begin{Eqn} \label{Tori}
T_x = \langle x,t \rangle; \ T_y = \langle y,t \rangle; \ T_z =
\langle z,t \rangle; \ T_w = \langle x=y=z, t \rangle.
\end{Eqn}

\noindent These tori are oriented as follows.  View $\T^4$ as $\T^3 \times
\sss^1$ and view the tori as given by taking products of circles in
$\T^3$ with the last factor.  The circles are projections of lines in
$\R^3$ which are all disjoint and are parallel to the $x,y,z$ axes and
to the line $x=y=z$ respectively.  Orient these lines by arrows
pointing to the positive infinity, and orient the circle
$\sss^1_{(t)}$; this orients the tori.  We choose the lines in $\R^3$
parallel to the obvious axes so that all the tori are disjoint.

\begin{Lem} \label{homologyrelation}
  There is a relation $[T_w] = [T_x + T_y + T_z]$ in $H_2 (\T^4;\Z)$.
\end{Lem}

\begin{Pf}
  This is straightforward from the given presentation.  We recall for
  the reader's convenience, however, the way this comes about via the
  original construction.  From \cite{McM-Taubes}: View $\T^3$ as
  surgery on the Borromean rings $\mathcal{B}$ in $\R^3$.  The circles
  $C_x, C_y, C_z$ in $\T^3$ for which $C_{\bullet} \times \sss^1 =
  T_{\bullet}$ are the images of meridians to the three zero-framed
  components of the link $\mathcal{B}$.  On the other hand, the circle
  $C_w$ is the image of an \emph{axis} to the link in $\R^3$, that is
  a line linking each component of the rings precisely once.  The
  meridians $C_x, C_y, C_z$ generate $H_1 (\T^3; \Z)$ and the
  co-ordinates of any loop $\gamma$ in $H_1$ with respect to this
  basis are given by the linking numbers, which yields the formula
  claimed.
\end{Pf}

\noindent From this construction, or directly, we can take two parallel and
unlinked copies of the axis in $\R^3$ and obtain two parallel copies
of the torus $T_w \subset \T^4$ each of which satisfies this homology
relation.  For the standard symplectic structure on $\T^4$ induced
from $dx \wedge dt + dy \wedge dz$ on $\R^4$, the tori $T_x$ and $T_w$
are symplectic whilst the tori $T_y, T_z$ are Lagrangian (but
canonically oriented by the orientations on the axes in $\R^3$).  By
the remark after (\ref{extend}), there is a uniquely defined oriented
smooth four-manifold $Z$ given by forming the fibre sum of $\T^4$ with
five copies of $E(1)$.  Here we glue a neighbourhood of a complex
oriented fibre $F \subset E(1)$ to the oriented tori $T_x, T_y, T_z$
and to two parallel copies of $T_w$.  Note that all the tori have
canonically trivial normal bundles, and this fixes the indeterminacy
in the fibre sum.  Because $\pi_1 (\T^4)$ is generated by the images
of the fundamental groups $\pi_1 (T_{\bullet}) \rightarrow \pi_1
(\T^4)$, and because $\pi_1 (E(1) \backslash F)$ is trivial (contract
any loop along a cusp fibre, which is a rational curve), it follows
easily that $Z$ is simply connected.  Note that the McMullen-Taubes
manifold of \cite{McM-Taubes} is exactly the smooth fibre sum of
$\T^4$ along the four oriented tori $T_x, T_y, T_z$ and one copy of
$T_w$.  Using our remarks on Chern classes, we have the following
proposition, from which the theorem (\ref{twoformsdiffer}) clearly
follows.

\begin{Prop}
  The homology classes $-3[T_x + T_y + T_z]$ and $-3[T_x + T_y] - T_z$
  are each the first Chern class of a symplectic structure on $Z$.
\end{Prop}

\begin{Pf}
  There is a symplectic structure $\omega_+$ on $\T^4$ for which all
  four tori are symplectic submanifolds with the symplectic
  orientation and given orientations agreeing.  Form the symplectic
  fibre sums with respect to such a form, inducing a form $\Omega_+$
  on $Z$.  Then we find by (\ref{choiceone}) that $c_1 (Z; \Omega_+)$
  is given by
  
  $$\sum_{i=1}^5 c_1 (\T^4; \omega_+) + [F_x + F_y + F_z + F_w + F_w]
  $$
  $$
  \qquad \qquad - 2 [F_x = T_x] - 2[F_y=T_y] - 2[F_z=T_z] -
  2[F_w=T_w] -2[F_w=T_w],$$

\noindent where we write $F_{\bullet}$ for the homology class of the fibre in
the $\bullet$-th copy of $E(1)$ for clarity.  However, we can also
choose a form $\omega_-$ on $\T^4$ for which $T_x, T_y$ and $T_w$ are
symplectic with the given orientation but for which $\omega_-$
restricts on $T_z$ to a form which is symplectic with the negative
orientation.  Now use (\ref{choicetwo}) to compute the Chern class of
the symplectic form $\Omega_-$ induced on $Z$ after performing the new
symplectic sum (which is smoothly identical to the old).  Recall from
(\ref{changesign}) that the sign of the first Chern class of $E(1)$
changes if we change the sign of the symplectic form: then $c_1 (Z;
\Omega_-)$ is now given by

$$\sum_{i=1}^5 c_1 (\T^4; \omega_-) + [F_x + F_y -F_z +F_w + F_w] $$
$$\qquad \qquad -2[F_x = T_x] - 2[F_y = T_y] + 2[F_z = T_z] - 2[F_w =
T_w] - 2[F_w=T_w].$$

\noindent Now recall that for both $\omega_{\pm}$ the first Chern class of
$\T^4$ vanishes, since these are perturbations of the standard form;
then expanding out and using (\ref{homologyrelation}) the proposition
follows.
\end{Pf}

\noindent One can fibre sum with further parallel copies of $T_w$ to
obtain a sequence of (pairwise non-homeomorphic)
simply connected manifolds with pairs of symplectic structures,
precisely one of which
has $c_1$ divisible by $n$, for any $n>2$.  If $n=p_1 \cdots p_N$ is a
product of primes and we perform $p_i$ of the 
fibre sums with the negatively oriented Lagrangian torus, we obtain a
symplectic structure with first Chern class divisible by $p_i$ but not
any other $p_{j \neq i}$.  This gives the main theorem
(\ref{alldimensions}).  For more variation,
replace the four-torus with some twisted torus bundle over the torus which
has a symplectic structure with vanishing first Chern class.  Infinite
families of such examples are described in \cite{ivangokova}.

%%%%%%%%%%%%%%%%%%%%%%%%%%%%%%%%%%%%%%%%%%%%%%%%%%%%%%%%%%%%%%%%%%%%%%

\bibliographystyle{amsplain}
\bibliography{main}

\end{document}